\documentclass{article}
\usepackage{emsnewsletter}

\usepackage[english]{babel}
\usepackage{hyperref}       
\usepackage{url}            
\usepackage{booktabs}       
\usepackage{amsfonts}       
\usepackage{tablefootnote}
\usepackage{verbatim}
\usepackage{nicefrac}       
\usepackage{microtype}      
\usepackage{lipsum}
\usepackage{mathtools}

\usepackage{amsmath,amssymb,amsthm}
\usepackage{algorithm,algorithmic}
\usepackage{pifont}
\usepackage{cases}
\usepackage{subcaption,graphicx}
\usepackage{stackengine}    
\usepackage{wrapfig}
\usepackage{enumitem}
\usepackage{multirow}
\usepackage{xcolor}
\usepackage{array}

\newtheorem{assumption}[theorem]{Assumption}
\newtheorem{defn}[theorem]{Definition}

\newcommand{\eps}{\varepsilon}

\newcommand{\ol}{\overline}
\newcommand{\one}{\mathbf{1}}

\renewcommand{\hat}{\widehat}

\DeclareMathOperator*{\argmin}{arg\,min}

\DeclareMathOperator{\spn}{span}
\DeclareMathOperator{\kernel}{Ker}
\DeclareMathOperator{\image}{Im}

\DeclareMathOperator{\col}{col}
\DeclareMathOperator{\gap}{Gap}

\newcommand{\R}{\mathbb{R}}
\newcommand{\Z}{\mathbb{Z}}

\newcommand{\E}{\mathbb{E}}

\newcommand{\bI}{{\bf I}}

\newcommand{\bL}{{\bf L}}

\newcommand{\bW}{{\bf W}}

\newcommand{\cE}{{\mathcal{E}}}

\newcommand{\cG}{{\mathcal{G}}}

\newcommand{\cL}{{\mathcal{L}}}

\newcommand{\cV}{{\mathcal{V}}}
\newcommand{\cW}{{\mathcal{W}}}

\newcommand{\cZ}{{\mathcal{Z}}}

\newcommand{\bx}{{\bf x}}
\newcommand{\by}{{\bf y}}
\newcommand{\bz}{{\bf z}}

\newcommand{\ds}{\displaystyle}
\newcommand{\norm}[1]{\left\| #1 \right\|}

\newcommand{\angles}[1]{\left\langle #1 \right\rangle}
\newcommand{\cbraces}[1]{\left( #1 \right)}
\newcommand{\sbraces}[1]{\left[ #1 \right]}
\newcommand{\braces}[1]{\left\{ #1 \right\}}
\def\<#1,#2>{\langle #1,#2\rangle}

\newcommand{\Lmax}{L_\ell}
\newcommand{\Lav}{L_g}
\newcommand{\Lf}{L_f}
\newcommand{\Lsummand}{L_s}
\newcommand{\Mmax}{M_\ell}
\newcommand{\Mav}{M_g}
\newcommand{\Mf}{M_f}
\newcommand{\mumin}{\mu_\ell}
\newcommand{\muav}{\mu_g}
\newcommand{\muf}{\mu_f}
\newcommand{\sigmamax}{\sigma_\ell}
\newcommand{\sigmaav}{\sigma_g}
\newcommand{\sigmamaxdual}{\sigma_{\ell,*}}
\newcommand{\sigmaavdual}{\sigma_{g,*}}

\newcommand{\blue}[1]{{\color{blue}#1}}

\newcommand{\rev}[1]{{\color{black}#1}}

\usepackage[colorinlistoftodos,bordercolor=blue,backgroundcolor=blue!20,linecolor=blue,textsize=scriptsize]{todonotes}

\author{
	\textit{Alexander Rogozin$^1$, Alexander Gasnikov$^{1,2,3}$, Aleksandr Beznosikov$^1$, Dmitry Kovalev$^4$}\\
	~1. Moscow Institute of Physics and Technology, Moscow, Russia \\
	~2. Caucasus Mathematical Center, Adyghe State University, Maikop, Russia\\
	~3. IITP RAS, Moscow, Russia\\
	~4. Universit\'e catholique de Louvain (UCL), CORE (IMMAQ), Louvain-la-Neuve, Belgium
}
\title{\Large{Decentralized Convex Optimization over Time-Varying Graphs: a Survey}}

\begin{document}
   
\maketitle

\begin{abstract}
	Decentralized optimization over time-varying networks has a wide range of applications in distributed learning, signal processing and various distributed control problems. The agents of the distributed system locally hold optimization objectives and can communicate to their immediate neighbors over a network that changes from time to time. \rev{In this paper, a survey of state-of-the-art results is presented and the techniques for optimization over time-varying graphs are described}. Moreover, an overview of open questions in the field is also given and hypotheses and directions for future work \rev{are formulated}.
\end{abstract}

\rev{

}

\section{Introduction}

Consider a sum-type optimization problem
\begin{align}\label{eq:sum_type_problem}
    \min_{x\in\R^d}~ f(x) = \frac{1}{m} \sum_{i=1}^m f_i(x).
\end{align}

The nodes can exchange information through a decentralized communication network, which may be time-varying. Each node is connected to several others via communication links and can communicate to them. A time-varying network is represented as a sequence of graphs (this paper focuses on undirected graphs). Since there is no centralized aggregator (server node, master node), each agent locally holds a copy $x_i$ of the decision vector $x$. The vectors held by the nodes should be synchronized, but the agents can communicate only to their immediate neighbors.

Problems of type \eqref{eq:sum_type_problem} arise in many applications where centralized aggregation is limited due to the structure of the network, privacy constraints or large amounts of data. Applications include vehicle coordination and control \cite{ren2008distributed}, distributed statistical inference and machine learning \cite{rabbat2004distributed,forero2010consensus,nedic2017fast}, power system control \cite{ram2009distributed,gan2012optimal}, distributed averaging \cite{cai2014average,olshevsky2014linear,xiao2007distributed}, formation control \cite{olshevsky2010efficient,ren2006consensus,jadbabaie2003coordination}, distributed spectrum sensing \cite{bazerque2009distributed}. See surveys \cite{nedic2010cooperative,nedic2020distributed} for additional examples \rev{and \cite{gorbunov2022recent,dvinskikh2019decentralized} for reviews of decentralized optimization over time-static graphs.}


Each $f_i$ is assumed to be convex and stored at a separate computational agent (or node). Moreover, \rev{each} $f_i$ is equipped with a first-order oracle (either stochastic or deterministic). That means that each computational node can compute gradients or stochastic gradients of the function it holds. During the computation process, the agents can exchange their decision vectors and gradients.

A decentralized optimization algorithm should be designed in such a way that the sum of functions is minimized while the decision vectors held at different computational nodes stay approximately the same. Assuming that node $i$ holds $x_i$, the optimal point in the decentralized sense should be consensual and optimal, i.e. 
$$
x_1 = \ldots = x_m = x^* = \argmin_{x\in\R^d} \frac{1}{m}\sum_{i=1}^m f_i(x).
$$


\subsection{Notation}

\rev{The column vector is denoted} $\col(x_1, \ldots, x_m) = (x_1^\top \ldots x_m^\top)^\top$. Let $x_i$ be held at node $i$ and denote the common decision vector $\bx = \col(x_1, \ldots, x_m)$. \rev{Also} introduce $\ol x = \frac{1}{m}\sum_{i=1}^m x_i$ and $\ol\bx = \col(\ol x, \ldots, \ol x)$. \rev{Moreover, denote} the identity matrix $\bI$ and the vector of ones $\one$ (the dimensions are known from the context) and let $\otimes$ denote the Kronecker product. \rev{Let $\spn(\cdot)$ denote} the linear span of vectors. Throughout the paper, $\norm{\cdot}$ denotes the Euclidean norm (the 2-norm).

The largest and minimal positive singular values of a matrix are denoted $\sigma_{\max}(\cdot)$ and $\sigma_{\min}^+(\cdot)$, respectively. The largest and minimal positive eigenvalues of a matrix are denoted $\lambda_{\max}(\cdot)$ and $\lambda_{\min}^+(\cdot)$, respectively. Euclidean ball of radius $r$ is denoted $B_r(x) = \{y\in\R^d:~ \norm{y - x}\leq r\}$.

\subsection{Solution accuracy}

\rev{Problem \eqref{eq:sum_type_problem} can be rewritten as} as
\begin{align}\label{eq:problem_consensus_constraints}
&\min_{\bx\in\R^{md}}~ F(\bx) = \sum_{i=1}^m f_i(x_i), \\
&\text{s.t.}~ x_1 = \ldots = x_m. \nonumber
\end{align}
Note that $\nabla F(\bx) = \col\cbraces{\nabla f_1(x_1) \ldots \nabla f_m(x_m)}$.

\rev{Vector $\hat\bx$ is called} an $\eps$-solution if 
\begin{align}\label{eq:solution_accuracy}
	&\frac{1}{m}\sum_{i=1}^m f_i(\hat x_i) - f^*\leq\eps,~~~ \cbraces{\frac{1}{m}\sum_{i=1}^m \norm{\hat x_i - \ol{\hat x}}^2}^{1/2}\leq C\eps.
\end{align}
Here $C$ is some constant not dependent on $\eps$. The constant $C$ may differ for different algorithms. Typically for dual methods constant $C$ depends on the norm of the dual problem solution. For each primal method, constant $C$ is individual.

For stochastic methods, the output $\hat\bx$ is stochastic, and by $\eps$-accuracy \rev{it is understood that}
\begin{align*}
	\E\sbraces{\frac{1}{m}\sum_{i=1}^m f_i(\hat x_i)} - f^*\leq \eps,~~~ \E\sbraces{\cbraces{\frac{1}{m}\sum_{i=1}^m \norm{\hat x_i - \ol{\hat x}}^2}^{1/2}}\leq C\eps.
\end{align*}

The complexity of distributed algorithms is measured by two quantities: the number of communication rounds and the number of oracle calls performed by a separate node.

\subsection{Paper organization}

In Section~\ref{sec:consensus}, \rev{the consensus over time-varying graphs is discussed and the assumptions on the communication network are formulated}. Section~\ref{sec:optimization} covers the main techniques and state-of-the-art-results in distributed optimization. After that, \rev{Section~\ref{sec:saddle} describes} the algorithms for decentralized saddle-point problems.

\section{Time-varying consensus}\label{sec:consensus}
The nodes communicate through a time-varying communication network. This paper is focused on \textit{undirected} graphs. The network is represented as a sequence of graphs $\{\cG^k = (\cV, \cE^k)\}_{k=0}^\infty$. \rev{Sequence $\{\cG^k\}_{k=0}^\infty$ is also referred to} as a \textit{time-varying graph}. The problems of reaching a consensus over time-varying graphs has been studied since 1980's (see i.e. seminal works \cite{tsitsiklis1984problems} and \cite{bertsekas1989parallel}). \rev{More recent works include \cite{ren2011distributed,ren2008distributed,proskurnikov2021delay,olfati2004consensus,kia2019tutorial}.}

In order to maintain consensus constraints $x_1 = \ldots = x_m$ in problem~\eqref{eq:problem_consensus_constraints}, a sequence of communication matrices is assigned to the time-varying graph. \rev{The two }most used types of communication matrices are \textit{mixing matrices} and \textit{gossip matrices}. Typically, mixing matrices are utilized in methods that use primal oracle and gossip matrices are employed in dual algorithms. Both types of matrices are defined in such way that a matrix-vector multiplication corresponds to one synchronized communication round.

\subsection{Mixing matrix}

\begin{assumption}\label{assum:mixing_matrix_sequence}
	Mixing matrix sequence $\{W^k\}_{k=0}^\infty$ satisfies the following properties.
	\item (Decentralized property) $[W^k]_{ij} = 0$ if $(i, j)\notin \cE^k$ and $i\ne j$.
	\item (Double stochasticity) $W^k \one = \one$ and $\one^\top W^k =\one^\top$.
	\item (Spectrum property) There exists a positive $\tau\in\Z$ and $\chi > 0$ such that for any $k\geq \tau - 1$ and any $x\in\R^d$ \rev{it holds}
	\begin{align*}
	\norm{\cbraces{W_\tau^k - \frac{1}{m}\one\one^\top}x}^2\leq \cbraces{1 - \frac{\tau}{\chi}}\norm{x}^2,
	\end{align*}
	where $W_\tau^k = W^k\ldots W^{k - \tau + 1}$.
\end{assumption}
\begin{remark}
	\rev{The} definition of $\chi$ \rev{given in Assumption~\ref{assum:mixing_matrix_sequence}} is slightly different from the one typically used in the literature. For example, \cite{nedic2017achieving} and \cite{kovalev2021lower} assume
	\begin{align*}
	\norm{\cbraces{W_\tau^k - \frac{1}{m}\one\one^\top}x}^2\leq \cbraces{1 - \frac{1}{\chi}}\norm{x}^2.
	\end{align*}
	\rev{Such} definition (i.e. \rev{$\chi$ can be called a} \textit{reduced condition number} of the time-varying network) \rev{is needed} solely for simplicity of notation: it allows to write $\chi$ instead of $\tau\chi$ in all complexity bounds. In the particular case of $\tau = 1$ (i.e. when the graph stays connected at every iteration), the definitions coincide.
\end{remark}
Also for each $k$ \rev{introduce} $\bW^k = W^k\otimes \bI$. The spectrum property in Assumption~\ref{assum:mixing_matrix_sequence} ensures geometric convergence of consensus iterates $\bx^{k+1} = \bW^k\bx^k$ to consensual point $\ol\bx^0$. After $N = O\cbraces{\chi\log\cbraces{\frac{1}{\eps}}}$ iterations \rev{it holds} $\norm{\bx^N - \ol\bx^0}^2\leq\eps$. Here and below the dependence from different parameters except $\eps$ under $\log\left(~\cdot~\right)$ \rev{is skipped}.

\begin{remark}[Sufficient conditions for Assumption~\ref{assum:mixing_matrix_sequence}]
	Paper \cite{nedic2017achieving} gives sufficient conditions for Assumption~\ref{assum:mixing_matrix_sequence} to hold. Firstly, the graph sequence $\{\cG^k\}_{k=0}^\infty$ should be $\tau$-connected. That means that for any $k\geq 0$ the union of $\tau$ consequent graphs $\hat\cG^k = \{\cV, \bigcup_{i=k}^{k + \tau - 1} \cE^{i}\}$ is connected. Secondly, the following restrictions on the mixing matrix weights are imposed:
	\item 1. (Double stochasticity) $W^k\one = \one,~ \one^\top W^k = \one^\top$.
	\item 2. (Positive diagonal) $[W^k]_{ii} > 0$ for $i = 1, \ldots, m$.
	\item 3. (Edge utilization) If $(i, j)\in\cE^k$, then $[W^k]_{ij} > 0$, else $[W^k]_{ij} = 0$.
	\item 4. (Nonvanishing weights) There exists $\theta> 0$ such that if $[W^k]_{ij} > 0$, then $[W^k]_{ij}\geq\theta$.
\end{remark}

In other words, the term $\tau$ in the Spectrum property of Assumption~\ref{assum:mixing_matrix_sequence} describes the number of iterations such that the union of $\tau$ consequent graphs is connected.

Typically mixing matrices are used in primal algorithms.

A possible way to build mixing matrices satisfying Assumption~\ref{assum:mixing_matrix_sequence} is to use Metropolis weights \cite{nedic2017achieving}:
\begin{align*}
	[W^k]_{ij} = 
	\begin{cases}
	\frac{1}{\max(\deg(i), \deg(j)) + 1} &\text{if } (i, j)\in \cE^k, \\
	0 &\text{if } (i, j)\notin \cE^k \text{ and } i\neq j, \\
	1 - \sum_{j\neq i} [W^k]_{ij}, &\text{$i = j$} \\	
	\end{cases}
\end{align*}

\subsection{Gossip matrix}

Dual methods typically use a notation of a \textit{gossip matrix} sequence $\{\cL^k\}_{k=0}^\infty$.

\begin{assumption}\label{assum:gossip_matrix_sequence}
	Gossip matrix sequence $\{\cL^k\}_{k=0}^\infty$ satisfies the following properties.
	\item 1. $[\cL^k]_{ij} = 0$ if $i\neq j$ and $(i, j)\notin\cE^k$.
	\item 2. $\kernel \cL^k \supseteq \spn(\one)$.
	\item 3. $\image \cL^k \subseteq \{x\in\R^m:~ x_1 + \ldots + x_m = 0\}$.
	\item 4. There exists a positive $\tau\in\Z$ and $\chi > 0$ such that for any $k\geq\tau - 1$ it holds 
	\begin{align*}
		\norm{\cL_\tau^kx - x}^2\leq \cbraces{1 - \frac{\tau}{\chi}}\norm{x}^2
	\end{align*} for all $x\in\R^m$ such that $x_1 + \ldots + x_m = 0$. Here $\cL_\tau^k$ is defined as
	\begin{align*}
		\bI - \cL_\tau^k = (\bI - \cL^k)\ldots (\bI - \cL^{k-\tau+1}).
	\end{align*}
\end{assumption}

For time-static networks, let $\cL$ denote the gossip matrix. Consensus constraints in problem \eqref{eq:problem_consensus_constraints} can be written as $\cL\bx = 0$. Therefore, distributed optimization is formulated as an affinely constrained problem. The dual approach build upon optimization of the function dual to $F$ s.t. $\cL\bx = 0$.

A possible way to obtain a gossip matrix is $\cL^k = \bI - W^k$. As noted in \cite{kovalev2021lower}, mixing matrix sequence $\{W^k\}_{k=0}^\infty$ satisfies Assumption~\ref{assum:mixing_matrix_sequence} if and only if gossip matrix sequence $\{\cL^k = \bI - W^k\}_{k=0}^\infty$ satisfies Assumption~\ref{assum:gossip_matrix_sequence}.

Alternatively, a gossip matrix can be built using a graph Laplacian. \rev{Matrix $\bL(\cG^k)$ is called a Laplacian of $\cG^k$} if
\begin{align*}
	\mathbf{L}(\cG^k) = 
	\begin{cases}
		\deg(i), &i = j, \\
		-1, &(i, j)\in\cE^k, \\
		0, &(i, j)\notin\cE^k \text{ and } i\ne j.
	\end{cases}
\end{align*}
Then $\braces{\cL^k = \bL(\cG^k)/\lambda_{\max}(\bL^k)}_{k=0}^\infty$ satisfies Assumption~\ref{assum:gossip_matrix_sequence}. Moreover, in the case $\tau = 1$ \rev{one} can equivalently define $\chi$ as
\begin{align*}
	\chi = \sup_{k\geq 0}\frac{\lambda_{\max}(\bL(\cG^k))}{\lambda_{\min}^+(\bL(\cG^k))}.
\end{align*}

\subsection{Multi-step consensus}

\rev{In this section} the multi-step consensus procedure \rev{is discussed}. Consider a problem with $L$-smooth and $\mu$-strongly convex summands $f_i$ distributed over network satisfying Assumption~\ref{assum:mixing_matrix_sequence} with parameters $\tau, \chi$. Let every iteration of the method corresponds to $O(1)$ local computations and $O(T)$ communications. To reach $\eps$-accuracy one needs to perform {\footnotesize$O\cbraces{\chi(T)/\tau\sqrt{L/\mu}\log(1/\eps)}$} local computations and {\footnotesize$O\cbraces{\chi(T)\sqrt{L/\mu}\log(1/\eps)}$} communications. However, the lower bound on the number of oracle calls is $O\cbraces{\sqrt{L/\mu}\log(1/\eps)}$. In order to control communication and computation complexities, a multi-step technique is used.

On the one hand, \rev{one} can choose $T = T_1 = \tau$. In this case $\chi(T_1) = \chi$ and communication complexity \rev{is} {\small$O\cbraces{\chi\sqrt{L/\mu}\log(1/\eps)}$} and local computation complexity $O\cbraces{\chi/\tau\sqrt{L/\mu}\log(1/\eps)}$. In the case of $\tau = 1$ the complexities coincide, which corresponds to a \textit{single-step consensus}, i.e. performing one communication step after each computation.

On the other hand, setting $T = T_2 = \lceil\chi\log 2\rceil$ \rev{leads to} $\chi(T_2) = O(1)$. Therefore, the required number of communications is $O\cbraces{\chi\sqrt{L/\mu}\log(1/\eps)}$ and the number of oracle calls is $O\cbraces{\sqrt{L/\mu}\log(1/\eps)}$.

If $\tau\geq 2$, both $T_1 = \tau$ and $T_2 = \lceil\chi\log 2\rceil$ lead to making several consensus iterations after each oracle call. However, in the case $\tau = 1$ \rev{note that} $T = T_1$ refers to a \textit{single-step consensus} while $T = T_2$ describes to a \textit{multi-step consensus}. The case $\tau = 1$ was historically the first to be studied in the literature since it particularly covers time-static graphs. Therefore, \rev{the case $T = T_1$ is referred to as a \textit{single-step consensus}} and \rev{the scheme with $T = T_2$ is referred to as} a \textit{multi-step consensus} even when $\tau\geq 2$.

Multi-step communication procedure is used in distributed optimization over static graphs, as well. In the time-static case such technique is called \textit{Chebyshev acceleration} and it allows to additionally reduce communication complexity to $\sqrt\chi$. Given a gossip matrix $W$ with condition number $\chi$, \rev{one can} replace $W$ by a Chebyshev polynomial $P_K(W)$ of degree $K = \lceil\chi\rceil$. The condition number of $P_K(W)$ equals $O(1)$ due to the specific structure of the polynomial (see \cite{scaman2017optimal}).

The acceleration on communication steps is not possible in the time-varying scenario, i.e. \rev{it is not possible to reduce} $\chi$ to $\sqrt\chi$ in the time-varying setting. This follows from lower complexity bounds \cite{kovalev2021lower}.

Reaching the consensus over time-varying networks can be viewed as a quadratic optimization problem with a time-varying objective function $\{c^k(\bx) = \frac{1}{2}\bx^\top(\bI - \bW^k)\bx\}_{k=0}^\infty$. All functions $c^k(\bx)$ have the same minimizer $\ol\bx^0$. Non-accelerated gradient descent corresponds to a consensus algorithm. On each iteration of gradient descent, there is a contraction of potential function $\Phi^k = \frac{1}{2}\norm{\bx^k - \ol\bx^0}^2$. This contraction is robust to graph changes, which makes non-accelerated consensus converge over time-varying graphs.

On the contrary, accelerated gradient methods build upon a potential function of type $\Phi^k = a_k(c(\by^k) - c^*) + b_k\norm{\bz^k - \ol\bx^0}^2$ \cite{bansal2019potential}, where $\by^k, \bz^k$ are additional extrapolation sequences and $a_k, b_k > 0$ are scalars (note that in \rev{the case of consensus problems} optimal value $c^* = 0$). The first summand in $\Phi^k$ contains the function value and therefore is not robust to network changes.

\section{Optimization}\label{sec:optimization}

In this section, several classes of objective functions \rev{are covered}. The results for smooth and non-smooth (strongly) convex functions with deterministic and stochastic gradients \rev{are described}. Both algorithms that use primal and dual oracle are covered. Moreover, the possible extensions to novel classes of problems \rev{such as} data similarity, directed graphs and alternative assumptions on the time-varying networks \rev{are also briefly discussed}.

\subsection{Definitions and assumptions}
Recall several standard definitions first.
\begin{defn}
	Consider function $h:~ \R^d\to\R$.
	\item 1. Function $h$ is convex if for any $x_1, x_2\in \R^d$ it holds
	\begin{align*}
		h(x_2)\geq h(x_1) + \angles{\nabla h(x_1), x_2 - x_1}.
	\end{align*}
	\item 2. Function $h$ is $\mu$-strongly convex if for any $x_1, x_2\in \R^d$ it holds
	\begin{align*}
	h(x_2)\geq h(x_1) + \angles{\nabla h(x_1), x_2 - x_1} + \frac{\mu}{2}\norm{x_2 - x_1}^2.
	\end{align*}
	\item 3. Function $h$ is $L$-smooth if for any $x_1, x_2\in \R^d$ it holds
	\begin{align*}
	h(x_2)\leq h(x_1) + \angles{\nabla h(x_1), x_2 - x_1} + \frac{L}{2}\norm{x_2 - x_1}^2.
	\end{align*}
\end{defn}

The assumptions on objective functions are standard for convex optimization.
\begin{assumption}\label{assum:convex}
	For every $i = 1, \ldots, m$, function $f_i$ is convex.
\end{assumption}

\begin{assumption}\label{assum:str_convex}
	For every $i = 1, \ldots, m$, function $f_i$ is $\mu_i$-strongly convex.
\end{assumption}

\begin{assumption}\label{assum:smooth}
	For every $i = 1, \ldots, m$, function $f_i$ is $L_i$-smooth.
\end{assumption}

\begin{assumption}\label{assum:bounded_grad}
	For every $i = 1, \ldots, m$, the norm of subgradients of $f_i$ is bounded by $M_i$, i.e.
	\begin{align*}
	\norm{\partial f_i(x)}\leq M_i.
	\end{align*}
\end{assumption}
Further in the paper, Assumptions~\ref{assum:convex}, \ref{assum:str_convex}, \ref{assum:smooth} and \ref{assum:bounded_grad} \rev{are used} in different combinations. The stochastic case \rev{is also covered}, for which the following assumption \rev{is needed}.

\begin{assumption}\label{assum:stoch_gradient}
	For every $i = 1, \ldots, m$, function $f_i$ is equipped with a non-biased stochastic gradient $\nabla f_i(x, \xi_i)$ with variance bounded by $\sigma_i^2$. For any $x\in\R^d$ it holds
	\begin{align*}
	\E\nabla f_i(x, \xi_i) = \nabla f_i(x),~~~ \E\norm{\nabla f_i(x, \xi_i) - \nabla f_i(x)}^2\leq \sigma_i^2.
	\end{align*}
	Here random variables $\{\xi_i\}_{i=1}^m$ are independent.
\end{assumption}

The results in the literature use local, global and original constants characterizing objective functions.

The worst-case constants of $f_i$ \rev{are referred to as \textit{local} constants}. That is, 
\begin{subequations}\label{eq:def_max_constants}
\begin{align}
\Lmax &= \max_{i = 1, \ldots, m} L_i,~ \Mmax = \max_{i = 1, \ldots, m} M_i,\\
\mumin &= \min_{i = 1, \ldots, m} \mu_i,~ \ds\sigmamax = \max_{i = 1, \ldots, m}\sigma_i^2.
\end{align}
\end{subequations}
Note that $L_\ell$ and $\mu_\ell$ are the smoothness and strong convexity constants of function $F$ defined in \eqref{eq:problem_consensus_constraints}, respectively.

\noindent \rev{\textit{Global} constants are understood as} the average constants, i.e.
\begin{subequations}\label{eq:def_av_constants}
\begin{align}
\Lav &= \frac{1}{m}\sum_{i=1}^m L_i,~ \Mav = \frac{1}{m}\sum_{i=1}^m M_i, \\
\muav &= \frac{1}{m}\sum_{i=1}^m \mu_i,~ \sigmaav^2 = \frac{1}{m}\sum_{i=1}^m \sigma_i^2.
\end{align}
\end{subequations}


\noindent \rev{\textit{Original} constants denote} the constants $\Lf, \muf, \Mf$ corresponding to function $f$ itself.

Note that $\Lf\leq\Lav\leq\Lmax,~ \Mf\leq\Mav\leq\Mmax,~ \muf\geq\muav\geq\mumin,~ \sigmaav^2\leq\sigmamax^2$.
\begin{remark}
	Local and global may significantly differ. For example, consider $f_i(x) = \norm{x}^2/2$ for $i = 1, \ldots, m - 1$ and $f_m(x) = m\norm{x}^2 / 2$. Then $\Lmax = m,~ \Lav = 2 - 1/m$. If the number of nodes is big, \rev{then} $\Lmax\gg\Lav$. Moreover, if $x$ is restricted to a Euclidean ball $B_{r}(0)$ for some $r > 0$, \rev{it also holds} $\Mmax = m\gg (2-1/m) = \Mav$.
	
	Note that global and original constants may differ as well. Let $m = d$, and $f_i(x) = (x^{(i)})^2/2$, where $x^{(i)}$ is the $i$-th component of $x$. Then $\Lav = 1,~ \Lf = 1/d$, i.e. $\Lav / \Lf = d$.
	
	Moreover, a trick suggested in \cite{scaman2017optimal} can improve the dependence from $\mumin$ to $\muav$. Namely, \rev{each} $f_i(x)$ \rev{can be replaced} by $\hat f_i(x) = f_i(x) + (\muav - \mu_i)/2 \norm{x}^2$.
\end{remark}

Also throughout this section, let $R = \norm{x^0 - x^*}$ denote the distance from initial guess to solution.
%

\subsection{Algorithm techniques}\label{subsec:techniques}

First, \rev{it is convenient to} show the concepts of decentralized optimization on a simple algorithm. Decentralized gradient descent~\cite{yuan2016convergence,nedic2009distributed,nedic2009subgradient} is a method that directly combines gradient steps and communication rounds.
\begin{align}\label{eq:dgd}
	\bx^{k+1} = \bW^k \bx^k - \gamma_k \nabla F(\bx^k).
\end{align}
According to decentralized property in Assumption~\ref{assum:mixing_matrix_sequence}, the $i$-th node computes $x_i^{k+1}$ only by communications with its immediate neighbors at time step $k$. Therefore, $\bW^k\bx^k$ corresponds to one (synchronous) communication round. Schemes of type \eqref{eq:dgd} are relatively simple to analyze (see i.e. analysis via Lyapunov functions in \cite{yuan2016convergence}) but do not reach optimal complexity bounds.

\rev{It is possible to enumerate the three techniques in the literature that led} to different state-of-the-art primal and dual algorithms: inexact oracle, gradient tracking and ADOM.

\textbf{Technique 1: gradient tracking}. A group of algorithms uses a \textit{gradient tracking} technique. This approach assumes that the local gradients at the nodes are averaged as well as decision vectors. For example, \rev{initialize $\bx^0$, put $\by^0 = \nabla F(\bx^0)$} and run procedure \cite{nedic2017achieving}
\begin{subequations}\label{eq:diging}
\begin{align}
	\bx^{k+1} &= \bW^k\bx^k - \gamma\by^k, \\
	\by^{k+1} &= \bW^k\by^k + \nabla F(\bx^{k+1}) - \nabla F(\bx^k).
\end{align}
\end{subequations}
Here $\by^k$ stands for gradient approximation. Gradient tracking allows each node to have a local approximation $y_i^k\approx (1/m)\sum_{i=1}^m \nabla f_i(x_i^k)$ of the average gradient over the network. The technique was successfully applied to time-varying networks in \cite{nedic2017achieving} to show first geometric convergence rates for time-varying graphs. After that, gradient tracking allowed to reach lower complexity bounds by Acc-GT algorithm proposed in \cite{li2021accelerated}. The paper \cite{ye2020multi} also introduced Mudag that reached optimal bounds up to $\log(\Lmax/\mumin)$ factor. Gradient tracking algorithms have a non-complicated structure but their analysis is quite involved.

It is also possible to use different mixing matrices for vectors and gradients as in push-pull gradient methods \cite{pu2020push,nedich2022ab}. Push-pull methods are capable of working over directed graphs.

\textbf{Technique 2: ADOM}. \rev{In this survey, technique ADOM is named after the algorithm where it was originally used \cite{kovalev2021adom}}. The technique of ADOM \rev{is based on interpreting decentralized communication as a compression operator and applying error-feedback}. Initially ADOM was used for dual oracle in \cite{kovalev2021adom} and then was employed for primal oracle in \cite{kovalev2022optimal_1} \rev{the corresponding primal method is called ADOM+. Both} ADOM and ADOM+ reach lower complexity bounds. \rev{On the one hand}, the derivation of both algorithms is natural and logical. On the other hand, the algorithms have quite a complicated structure.

\textbf{Technique 3: inexact oracle}. One of the possible ways to develop a distributed algorithm is based on \textit{inexact oracle} concept \cite{devolder2014first}. After the nodes make a local computation step, the decision vectors held by the nodes are averaged up to target accuracy $\eps$ by running multi-step consensus for $T = O(\chi\log(1/\eps))$ iterations. According to Assumption~\ref{assum:mixing_matrix_sequence}, this allows to get a projection on the consensus constraint set $x_1 = \ldots = x_m$ with accuracy $\eps$. As a result, at each computation step the values at the nodes are approximately averaged, and therefore the gradient is obtained with approximation accuracy $\eps$. The inexact oracle concept allows to tackle the inexactness of the gradient.

Initially the inexact oracle technique was developed for time-static networks in \cite{jakovetic2014fast}. Its generalization on time-varying networks for accelerated gradient method were made in \cite{rogozin2021towards,dvinskikh2019decentralized,li2020decentralized}.

The main advantage of using the inexact oracle framework is the simplicity of interpretation and analysis of the distributed method. However, all such schemes have an additional logarithmic term in communication complexity (i.e. the complexity for strongly convex smooth objectives is proportional to $\log^2(1/\eps)$). Inexact oracle schemes need to perform a comparatively large number of communications between steps related to other techniques, and the number of communications needs to be accurately tuned.

\textbf{Regularization}. Some of the papers cover only the strongly-convex case. Such results can be generalized to a non-strongly convex problems using regularization.

\begin{lemma}[Regularization]\label{lem:regularization}\cite{gorbunov2019optimal}
	Let $h(x)$ be a convex function and $\eps > 0$ be the desired accuracy. Let $h(x)$ have a minimizer $x^*$ such that $\norm{x^*}\leq R$. Consider regularized function $\hat h(x) = h(x) + \frac{\eps}{2R^2}\norm{x}^2$ \rev{and denote its optimal value $\hat h^*$}. Let $\hat x$ be an $\eps/2$-minimizer of $\hat h$, i.e.
	\begin{align*}
		\hat h(\hat x) - \hat h^*\leq \frac{\eps}{2}.
	\end{align*}
	Then \rev{it holds} $h(\hat x) - h^*\leq \eps$.
\end{lemma}

In the following tables, the existing results \rev{are covered along with} their possible extensions. Some extensions can be maid by regularization, and some are just guesses based on distributed optimization over time-static graphs. \blue{\rev{The} hypotheses are marked in blue}.

\subsection{Primal oracle}

\rev{First, consider} primal algorithms and cover deterministic and stochastic cases. Lower bounds for problems on time-varying graphs satisfying Assumptions~\ref{assum:mixing_matrix_sequence}, \ref{assum:str_convex} and \ref{assum:smooth} are 
\begin{align*}
&\Omega\cbraces{\chi\sqrt{\frac{\Lav}{\muav}}\log\frac{1}{\eps}}~~ \text{communications}, \\
&\Omega\cbraces{\sqrt{\frac{\Lav}{\muav}}\log\frac{1}{\eps}}~~ \text{local oracle calls}
\end{align*}
as shown in \cite{kovalev2021lower}.

The constant $C$ used in \eqref{eq:solution_accuracy} to measure consensus is individual for every primal method. For example, for stochastic APM-C \cite{rogozin2021accelerated} \rev{it holds} $C = (m/32)~ \muav^{3/2}/(\Lav^{1/2}\Lmax^2)$.

The primal approach assumes that every node $i$ has access only to the gradient of (primal) function $f_i$.

The following table describes the results for primal oracle as well as \rev{the} hypotheses (\blue{in blue}). \rev{Notation $O\cbraces{\cdot}$ is omitted for brevity}.
\begin{table}[H]
	\begin{tabular}{|b{0.1cm}|b{0.9cm}|>{\centering\arraybackslash}b{2.9cm}|>{\centering\arraybackslash}b{2.9cm}|}
		\hline
		& & convex & str. convex \\ \hline
		\parbox[t]{20mm}{\multirow{2}{*}{\rotatebox[origin=r]{90}{\footnotesize smooth\hspace{0.5cm}}}} & & \blue{Acc-GT} \cite{li2021accelerated} & Acc-GT \cite{li2021accelerated} \\
		& & \blue{{ADOM+ \cite{kovalev2021lower}}} & {ADOM+ \cite{kovalev2021lower}} \\
		& ~comm. & \blue{$\chi\sqrt{\frac{\Lmax R^2}{\eps}}\log\frac{1}{\eps}$} & $\chi\sqrt{\frac{\Lmax}{\mumin}}\log\frac{1}{\eps}$ \\
		& ~comp. & \blue{$\sqrt{\frac{\Lmax R^2}{\eps}}\log\frac{1}{\eps}$} & $\sqrt{\frac{\Lmax}{\mumin}}\log\frac{1}{\eps}$ \\ \hline
		\vspace{0.2cm}
		\parbox[t]{20mm}{\multirow{2}{*}{\rotatebox[origin=c]{90}{\footnotesize{non-smooth}}}} & comm. & \blue{$\chi\frac{\Mav R}{\eps}\log\frac{1}{\eps}$} & \blue{$\chi\frac{\Mav}{\sqrt{\mumin\eps}}\log\frac{1}{\eps}$} \\ [0.25cm]
		& comp. & \blue{$\frac{\Mav^2R^2}{\eps^2}$} & \blue{$\frac{\Mav^2}{\mumin\eps}$} \\ [0.25cm] \hline
	\end{tabular}
	\caption{Primal oracle}
	\label{tab:optim_primal}
\end{table}

\textbf{Deterministic smooth objectives}. The bounds for smooth strongly convex scenario are achieved by ADOM+ \cite{kovalev2021lower} and Acc-GT \cite{li2021accelerated}. These algorithms appeared approximately at the same time and are based on different techniques. Acc-GT uses gradient tracking and Nesterov acceleration, while ADOM+ is based on the independent technique with error feedback and specific problem reformulation. The bounds for a smooth (non-strongly) convex case are obtained via regularization according to Lemma~\ref{lem:regularization}, which is maid in a straightforward way.

\textbf{Deterministic non-smooth objectives}. One of the first papers on distributed subgradient methods \cite{nedic2009distributed} proposed a method of type \eqref{eq:dgd} that converges to a $O(\gamma)$ neighborhood of the solution at a rate of $O(1/\eps)$, where $\gamma$ denotes the stepsize. These results were generalized to directed graphs in \cite{nedic2014distributed}. For the time-static scenario, optimal algorithms for non-smooth objectives \cite{dvinskikh2021decentralized} are derived using \textit{gradient sliding} technique \cite{lan2016gradient}. \rev{Guesses} are based on the results for time-static setup. The direct application of gradient sliding to time-varying networks is an open question. An extension to functions satisfying Polyak--Lojasiewicz condition was proposed in \cite{kuruzov2022gradient}.

Optimization of stochastic objectives over time-varying networks is not studied as good as with deterministic objectives. 

\textbf{Stochastic smooth objectives}. SGD with gradient tracking was studied in \cite{pu2021distributed}; however, it only converges to a neighborhood of the solution at suboptimal rates. Paper \cite{rogozin2021accelerated} proposed a penalty method (APM-C) that is based on inexact oracle framework and includes a batching technique. APM-C reaches lower bounds up to a $\log(1/\eps)$ factor. Its results can be generalized to convex case by using regularization (see Lemma~\ref{lem:regularization}).

\textbf{Stochastic non-smooth} objectives do not seem to be covered in the literature. \rev{The guesses are built} by analogy with time-static scenario and replacing $\sqrt{\chi}$ with $\chi$.

\begin{table}[H]
	\begin{tabular}{|b{0.1cm}|b{0.9cm}|>{\centering\arraybackslash}b{2.9cm}|>{\centering\arraybackslash}b{2.9cm}|}
		\hline
		& & convex & str. convex \\
		\hline
		\parbox[t]{20mm}{\multirow{2}{*}{\rotatebox[origin=r]{90}{\footnotesize smooth\hspace{0.3cm}}}} & & \blue{APM-C \cite{rogozin2021accelerated}} & APM-C \cite{rogozin2021accelerated} \\
		& ~comm. & \blue{$\chi\sqrt{\frac{\Lav R^2}{\eps}}\log^2\frac{1}{\eps}$} & $\chi\sqrt{\frac{\Lav}{\muav}}\log^2\frac{1}{\eps}$ \\
		& ~comp. & \blue{\footnotesize $\max\cbraces{\frac{\sigmaav^2 R^2}{m\eps^2}, \sqrt{\frac{\Lav R^2}{\eps}}\log\frac{1}{\eps}}$} & {\footnotesize $\max\cbraces{\frac{\sigmaav^2}{m\muav\eps}, \sqrt{\frac{\Lav}{\muav}}\log\frac{1}{\eps}}$} \\
		\hline
		\vspace{0.2cm}
		\parbox[t]{20mm}{\multirow{2}{*}{\rotatebox[origin=r]{90}{\footnotesize non-smooth\hspace{-0.15cm}}}} & ~comm. & \blue{$\chi\frac{\Mav R}{\eps}$} & \blue{$\chi\frac{\Mav}{\sqrt{\mumin\eps}}\log\frac{1}{\eps}$} \\ [0.3cm]
		& ~comp. & \blue{$\max\cbraces{\frac{\Mav^2 R^2}{\eps^2}, \frac{\sigmaav^2R^2}{m\eps^2}}$} & \blue{$\max\cbraces{\frac{\Mav^2}{\mumin\eps}, \frac{\sigmaav^2}{m\mumin\eps}}$} \\ [0.3cm]
		\hline
	\end{tabular}
	\caption{Primal stochastic oracle}
	\label{tab:optim_primal_stoch}
\end{table}

\subsection{Dual oracle}

Dual approach is an alternative to primal. It uses a different oracle, i.e. the gradient of the conjugate function.
\begin{definition}
	Consider function $h:~ \R^d\to\R$. \rev{Function} $h^*(y)$ is a conjugate of $h$ if
	\begin{align*}
	h^*(y) = \max_{x\in\R^d}\cbraces{\angles{x, y} - h(x)}.
	\end{align*}
\end{definition}
The dual approach in decentralized optimization utilizes a gossip matrix to interpret problem~\eqref{eq:problem_consensus_constraints} as a problem with affine constraints. The constraints are used to build a dual problem and solve it with gradient methods. Consider a time-static case with gossip matrix $\cL$. Dual problem to \eqref{eq:problem_consensus_constraints} writes as
\begin{align*}
	\min_{\by\in\R^{md}} \Phi(\by) = \max_{\bx\in\R^{md}} \sbraces{\angles{\by, \cL\bx} - F(\bx)}.
\end{align*}
One may see that computation of $\nabla\Phi(\by)$ corresponds to an oracle call and a communication round \cite{maros2018panda,scaman2017optimal,kovalev2021adom}.

For the dual approach over time-static networks, it can be shown that constant $C$ from \eqref{eq:solution_accuracy} used to bound consensus accuracy of the solution has form $C = 1 / \Mav$.

Typically, the following lemma is used to \rev{relate} optimization properties of the dual to the properties of the primal problem and therefore obtain convergence rates.
\begin{lemma}\label{lem:dual_smoothness_strong_convexity}\cite{rockafellar2015}
	Let function $h^*$ be a conjugate of $h$. If $h$ is $L$-smooth, then $h^*$ is $(1/L)$-strongly convex. If $h$ is $\mu$-strongly convex, then $h^*$ is $(1/\mu)$-smooth.
\end{lemma}

For stochastic objectives, \rev{it is assumed} that functions $f_i^*$ are equipped with a non-biased stochastic first-order oracle with bounded variance.
\begin{assumption}
	For every $i = 1, \ldots, m$, function $f_i^*$ is equipped with a non-biased stochastic gradient $\nabla f_i^*(y, \xi_i)$ with variance bounded by $\sigma_{i,*}^2$. For any $y\in\R^d$ it holds
	\begin{align*}
	\E\nabla f_i^*(y, \xi_i) = \nabla f_i^*(y),~~~ \E\norm{\nabla f_i^*(y, \xi_i) - \nabla f_i^*(y)}^2\leq \sigma_{i,*}^2.
	\end{align*}
	Here random variables $\{\xi_i\}_{i=1}^m$ are independent.
\end{assumption}
Analogously \rev{introduce} worst-case variance \\$\ds\sigmamaxdual^2 = \max_{i=1, \ldots, m} \sigma_{i,*}^2$ and average variance $\sigmaavdual^2 = (1/m)\sum_{i=1}^m \sigma_{i,*}^2$.

\textbf{Deterministic smooth objectives}. One of the first proposed linearly convergent dual methods for time-varying graphs is PANDA \cite{maros2018panda} that achieves a suboptimal communication complexity $O\cbraces{\chi^2\cbraces{\Lmax/\mumin}^{3/2}\log(1/\eps)}$ (using multi-step consensus will replace $\chi^2$ with $\chi$). PANDA applies gradient tracking technique to gradients of the dual function. An optimal dual algorithm is ADOM \cite{kovalev2021adom}. ADOM uses a specific reformulation of the problem. In particular, it treats multiplication by gossip matrix as a compression operator and employs an error feedback trick.

\textbf{Deterministic non-smooth objectives}. Non-accelerated methods for non-smooth problems with dual oracle were proposed in \cite{wu2019fenchel}. After that, ADOM was generalized to strongly convex non-smooth problems \cite{yufereva2022decentralized} by using Moreau--Yosida regularization. This generalization also covers constrained sets and supports non-Euclidean setup for the Wasserstein barycenter problem.

\begin{table}[H]
	\begin{tabular}{|b{0.1cm}|b{0.9cm}|>{\centering\arraybackslash}b{2.9cm}|>{\centering\arraybackslash}b{2.9cm}|}
		\hline
		& & convex & str. convex \\
		\hline
		\parbox[t]{20mm}{\multirow{2}{*}{\rotatebox[origin=r]{90}{\footnotesize smooth\hspace{0.3cm}}}} & & \blue{ADOM} \cite{kovalev2021adom} & ADOM \cite{kovalev2021adom} \\
		& ~comm. & \blue{$\chi\sqrt{\frac{\Lmax R^2}{\eps}}\log\frac{1}{\eps}$} & $\chi\sqrt{\frac{\Lmax}{\mumin}}\log\frac{1}{\eps}$ \\
		& ~comp. & \blue{$\sqrt{\frac{\Lmax R^2}{\eps}}\log\frac{1}{\eps}$} & $\sqrt{\frac{\Lmax}{\mumin}}\log\frac{1}{\eps}$ \\ \hline
		\parbox[t]{20mm}{\multirow{2}{*}{\rotatebox[origin=c]{90}{\footnotesize non-smooth\hspace{0.23cm}}}} & & \blue{ADOM \cite{yufereva2022decentralized}} & ADOM \cite{yufereva2022decentralized} \\
		& ~comm. & \blue{$\chi\frac{\Mav R}{\eps}\log\frac{1}{\eps}$} & $\chi\frac{\Mav}{\sqrt{\mumin\eps}}\log\frac{1}{\eps}$ \\ [0.2cm]
		& ~comp. & \blue{$\frac{\Mav R}{\eps}\log\frac{1}{\eps}$} & $\frac{\Mav}{\sqrt{\mumin\eps}}\log\frac{1}{\eps}$ \\ [0.2cm]
		\hline
	\end{tabular}
	\caption{Dual oracle}
	\label{tab:optim_dual}
\end{table}

Both ADOM and its generalization can be adapted for non-strongly-convex setup by using regularization (see \blue{hypotheses} in Table~\ref{tab:optim_dual}). Note that the analysis in \cite{yufereva2022decentralized} suggests that complexity bounds include $\Mmax$, but \rev{it may be possible to improve this dependence to $\Mav$}.

\textbf{Stochastic smooth and non-smooth objectives}. Even in the time-static scenario, \rev{it seems that only results for quadratic objectives \cite{zhang2022dual} are present in the literature}. Based on these results, \rev{the} \blue{hypotheses} \rev{are formulated} in Table~\ref{tab:optim_dual_stoch}.
\begin{table}[H]
	\begin{tabular}{|b{0.1cm}|b{0.9cm}|>{\centering\arraybackslash}b{2.9cm}|>{\centering\arraybackslash}b{2.9cm}|}
		\hline
		& & convex & str. convex \\
		\hline
		\parbox[t]{20mm}{\multirow{2}{*}{\rotatebox[origin=r]{90}{\footnotesize smooth\hspace{-0.05cm}}}} & ~comm. & {\footnotesize\blue{$\chi\sqrt{\frac{\Lmax R^2}{\eps}}$}} & {\footnotesize\blue{$\chi\sqrt{\frac{\Lmax}{\mumin}}\log\frac{1}{\eps}$}} \\
		& ~comp. & {\footnotesize \blue{$\chi\max\cbraces{\sqrt{\frac{\Lmax R^2}{\eps}}, \frac{\sigmaavdual^2}{m\mumin\eps}}$}} & {\footnotesize\blue{$\chi\max\cbraces{\sqrt{\frac{\Lmax}{\mumin}}\log\frac{1}{\eps}, \frac{\sigmaavdual^2}{m\mumin\eps}}$}} \\
		\hline
		\vspace{0.1cm}
		\parbox[t]{20mm}{\multirow{2}{*}{\rotatebox[origin=r]{90}{\footnotesize non-smooth\hspace{-0.2cm}}}} & ~comm. & {\footnotesize\blue{$\chi\frac{\Mav R}{\eps}$}} & {\footnotesize\blue{$\chi\frac{\Mav}{\sqrt{\mumin\eps}}\log\frac{1}{\eps}$}} \\ [0.2cm]
		& ~comp. & {\footnotesize \blue{$\chi\max\cbraces{\frac{\Mav R}{\eps}, \frac{\sigmaavdual^2\Mav^2}{m\eps}}$}} & {\footnotesize \blue{$\chi\max\cbraces{\frac{\Mav}{\sqrt{\mumin\eps}}, \frac{\sigmaavdual^2\Mav^2}{m\eps}}$}} \\ [0.2cm]
		\hline
	\end{tabular}
	\caption{Dual stochastic oracle}
	\label{tab:optim_dual_stoch}
\end{table}
Note that the bounds \cite{zhang2022dual} use worst-case variance $\sigmamaxdual^2$ but \rev{the guess is optimistic so} that it is possible to get a bound using average variance $\sigmaavdual^2$. 

Dual stochastic algorithms are an open direction for research.

\subsection{Finite-sum problems}\label{subsec:finite_sum_optimization}

Consider a finite-sum form of problems \eqref{eq:sum_type_problem}. Each $f_i(x)$ is now represented as a sum of $n$ functions.
\begin{align}\label{eq:sum_type_problem_finite_sum}
\min_{x\in\R^d}\max_{y\in\R^d}~ f(x) = \frac{1}{m} \sum_{i=1}^m \frac{1}{n} \sum_{j=1}^n f_{ij}(x).
\end{align}
Problems of type \eqref{eq:sum_type_problem_finite_sum} can be efficiently solved by stochastic methods that use additional knowledge of the sum-type structure. This class of algorithms is referred to as \textit{variance reduction} methods \cite{johnson2013accelerating,allen2016katyusha}. These methods compute the full gradient sum over all $f_i$ once in several iterations and make stochastic updates to the full gradient approximation afterwards.

The complexity of variance reduction methods typically depends on the worst-case smoothness constant over summands. Let $L_{ij}$ be the smoothness constant of $f_{ij}$ and introduce
\begin{align}\label{eq:def_Lsummand}
\Lsummand = \max_{\substack{i = 1, \ldots, m\\j = 1, \ldots, n}} L_{ij}.
\end{align}
In distributed optimization, variance reduction is done for time-static graphs \cite{hendrikx2020optimal}. The generalization to time-varying graphs is an open question. \rev{Probably} a method for decentralized finite-sum optimization of smooth functions over time-varying networks \rev{may be obtained} by combining Katyusha~\cite{allen2016katyusha} and ADOM/ADOM+ \cite{kovalev2021lower} or Acc-GT \cite{li2021accelerated}.
\begin{table}[H]
	\begin{tabular}{|b{1.0cm}|>{\centering\arraybackslash}b{3.1cm}|>{\centering\arraybackslash}b{3.1cm}|}
		\hline
		& convex-concave & str. convex-concave \\ \hline
		~comm. & \blue{$\chi\sqrt{\frac{\Lsummand R^2}{\eps}}\log\frac{1}{\eps}$} & \blue{$\chi\sqrt{\frac{\Lsummand}{\muav}}\log\frac{1}{\eps}$} \\
		~comp. & \blue{$\sqrt{n}\sqrt{\frac{\Lsummand R^2}{\eps}}\log\frac{1}{\eps}$} & \blue{$\sqrt{n}\sqrt{\frac{\Lsummand}{\muav}}\log\frac{1}{\eps}$} \\
		\hline
	\end{tabular}
	\caption{Finite-sum smooth optimization problems}
	\label{tab:finite_sum_optimization}
\end{table}

\subsection{Time-varying network as a stochastic sequence}

Most of the results of this paper are derived in the deterministic scenario: that is, a mixing (or gossip) matrix sequence $\{W^k\}_{k=0}^\infty$ is deterministic and satisfies Assumption~\ref{assum:mixing_matrix_sequence}. Alternatively, each of the mixing matrices $W^k$ may be seen as a random sample from some distribution $\cW^k$. As suggested in \cite{koloskova2020unified}, \rev{it can be assumed that} for all $\bx\in\R^{md}$ it holds
\begin{align}\label{eq:mixing_matrix_sequence_stochastic}
	\E\sbraces{\norm{W_\tau^k\bx - \ol\bx}^2}\leq \cbraces{1 - \frac{\tau}{\chi}}\norm{\bx - \ol\bx}^2,
\end{align}
where the expectation is taken over the distributions of $W^t$ and indices $t = k, \ldots, k - \tau + 1$. Such assumption broadens the class of networks compared to Assumption~\ref{assum:mixing_matrix_sequence}. In particular, it admits a scheme when only two randomly chosen agents communicate at once (see \cite{boyd2006randomized} for details) and also includes Local SGD schemes \cite{stich2018local,woodworth2020local,lin2018don}. Distributed SGD under Assumption~\eqref{eq:mixing_matrix_sequence_stochastic} was analyzed in \cite{koloskova2020unified} for optimization and in \cite{beznosikov2021decentralized} for saddle-point problems. Also \cite{trimbach2021acceleration} applied Catalyst to the results of \cite{koloskova2020unified} to get acceleration.

\subsection{Statistical similarity}

The method of using convex surrogate functions at local computation was proposed in \cite{sun2016distributed}. Combined with gradient tracking technique, it yields a SONATA algorithm \cite{sun2016distributed,scutari2019distributed} that supports non-convex objectives and works over time-varying \textit{directed} graphs. Several variants of SONATA \cite{sun2022distributed,scutari2019distributed,sun2016distributed} also support constrained optimization.

In the undirected time-static case with convex objectives, SONATA exploits statistical similarity of local functions \cite{sun2022distributed,tian2022acceleration}. That is, for any $x$ and $i = 1, \ldots, m$ \rev{one has}
\begin{align*}
\norm{\nabla^2 f_i(x) - \nabla^2 f(x)}\leq \beta
\end{align*}
for some $\beta > 0$. Under similarity assumption, SONATA communication rates (in time-static setup) are improved to $O(\sqrt{\chi}\beta/\muf\log(1/\eps))$ \cite{sun2022distributed}. Additionally, with Catalyst acceleration, the rates are further improved to $O\cbraces{\sqrt{\chi}\sqrt{\beta/\muf}\log^2(1/\eps)}$ \cite{tian2022acceleration}. Optimal algorithms were developed in \cite{kovalev2022optimal_2}. \rev{Note} that exploiting statistical similarity in the time-varying case is an open question. Using SONATA and the technique of convex surrogate functions may be a good starting point in this direction.

\subsection{Problems with affine constraints}

Affinely constrained problems over time-varying networks are still another direction of research. For problems of type
\begin{align*}
	&\min_{x\in\R^d}~ f(x) = \frac{1}{m}\sum_{i=1}^m f_i(x), \\
	&\text{s.t. } Ax = b
\end{align*}
in the time-static case a range of algorithms was proposed in \cite{necoara2011parallel,necoara2014distributed,necoara2015linear,yarmoshik2022decentralized,rogozin2022decentralized}. The algorithms are based either on dual decomposition or on primal approaches to affinely constrained minimization \cite{kovalev2021accelerated}. Generalization of the results to time-varying architectures is an open branch of research.

%
%
%

\section{Saddle-point problems}\label{sec:saddle}

In this section, sum-type min-max problems \rev{are considered}.
\begin{align}\label{eq:sum_type_saddle}
	\min_{x\in\R^{d_x}} \max_{y\in\R^{d_y}}~ f(x, y) = \frac{1}{m} \sum_{i=1}^m f_i(x, y).
\end{align}

\subsection{Definitions and assumptions}
\rev{Recall} the standard definitions for saddle-point problems. 

\begin{defn}
	Consider function $h(x, y)$. Introduce $z = \col(x, y)$ and a vector field $g(z) = \col\cbraces{\nabla_x h(x, y), -\nabla_y h(x, y)}$ associated with function $h(x, y)$.
	\item 1. Function $h(x, y)$ is convex-concave if the corresponding vector field $g$ is monotone, i.e. for all $z_1, z_2\in\R^{d_x + d_y}$ it holds
	\begin{align*}
		\angles{g(z_2) - g(z_1), z_2 - z_1}\geq 0.
	\end{align*}
	\item 2. Function $h(x, y)$ is $\mu$-strongly convex-concave if the corresponding vector field $g$ is $\mu$-strongly monotone, i.e. for all $z_1, z_2\in\R^{d_x + d_y}$ it holds
	\begin{align*}
	\angles{g(z_2) - g(z_1), z_2 - z_1}\geq \mu\norm{z_2 - z_1}^2.
	\end{align*}
	\item 3. Function $h(x, y)$ is $L$-smooth if the corresponding vector field $g$ is $L$-Lipschitz, i.e. for all $z_1, z_2\in\R^{d_x + d_y}$ it holds
	\begin{align*}
	\norm{g(z_2) - g(z_1)}\leq L\norm{z_2 - z_1}.
	\end{align*}
\end{defn}

\begin{assumption}\label{assum:convex_concave}
	For each $i = 1, \ldots, m$ function $f_i(x, y)$ is convex-concave.
\end{assumption}

\begin{assumption}\label{assum:stongly_convex_stongly_concave}
	For each $i = 1, \ldots, m$ function $f_i(x, y)$ is $\mu_i$-strongly convex-concave.
\end{assumption}

\begin{assumption}\label{assum:saddle_smooth}
	For each $i = 1, \ldots, m$ function $f_i(x, y)$ is $L_i$-smooth.
\end{assumption}

\begin{assumption}\label{assum:bounded_gradient}
	For each $i = 1, \ldots, m$ function $f_i(x, y)$ has a bounded gradient, i.e. there exists $M_i > 0$ such that
	\begin{align*}
		\norm{\nabla_x f_i(x, y)}\leq M_i,~~\norm{\nabla_y f_i(x, y)}\leq M_i.
	\end{align*}
\end{assumption}

Similarly to Section~\ref{sec:optimization}, \rev{studying} stochastic min-max problems \rev{requires} the following assumption.
\begin{assumption}\label{assum:saddle_stoch_gradient}
	For each $i = 1, \ldots, m$ function $f_i(x, y)$ is equipped with a non-biased gradient oracle with bounded variance, i.e. there exists $\sigma_i^2 > 0$ such that
	{\small
	\begin{align*}
		&\E\nabla_x f_i(x, y, \xi_i) = \nabla_x f_i(x, y),~~ \E\nabla_y f_i(x, y, \xi_i) = \nabla_y f_i(x, y), \\
		&\E\cbraces{\norm{\nabla_x f_i(x, y, \xi_i) - \nabla f_i(x, y)}^2 + \norm{\nabla_y f_i(x, y, \xi_i) - \nabla_y f_i(x, y, \xi_i)}^2}\leq \sigma_i^2.
	\end{align*}
	}
	Random variables $\{\xi_i\}_{i=1}^m$ are independent.
\end{assumption}
Similarly to Section~\ref{sec:optimization}, introduce worst-case constants $\Lmax, \Mmax$, $\mumin, \sigmamax$ as in \eqref{eq:def_max_constants}, average constants $\Lav, \Mav, \muav, \sigmaav$ as in \eqref{eq:def_av_constants} and constants that describe $f$ itself $\Lf, \Mf, \muf$.

For saddle-point problems, introduce a measure of solution accuracy.
\begin{defn}
	Consider function $h(x, y)$, vector field $g(z)$ associated with it and set $\cZ\subseteq\R^{d_x+d_y}$. For a point $\hat z = \col(\hat x, \hat y)\in \R^{d_x + d_y}$, introduce dual gap as follows.
	\begin{align*}
		\gap_\cZ^h(\hat z) = \max_{z\in\cZ}~ \angles{g(z), z - \hat z}.
	\end{align*}
\end{defn}
In order to measure solution accuracy, fix $D > 0$ and \rev{let} $(\hat\bx, \hat\by)$ \rev{be called} an $\eps$-solution of a distributed saddle-point problem if
\begin{align*}
	\gap^f_{B_D(0)}\cbraces{\ol{\hat x}, \ol{\hat y}}\leq\eps,~~ \cbraces{\frac{1}{m}\sum_{i=1}^m \norm{\hat x_i - \ol{\hat x}}^2 + \norm{\hat y_i - \ol{\hat y}}^2}^{1/2}\leq C\eps.
\end{align*}
Similarly to optimization methods, constant $C$ does not depend on $\eps$ but may be individual for every algorithm. If the output of the method is stochastic, \rev{by $\eps$-solution it is meant}
{\small
\begin{align*}
\E\sbraces{\gap^f_{B_D(0)}\cbraces{\ol{\hat x}, \ol{\hat y}}}\leq\eps,~~ \E\sbraces{\cbraces{\frac{1}{m}\sum_{i=1}^m \norm{\hat x_i - \ol{\hat x}}^2 + \norm{\hat y_i - \ol{\hat y}}^2}^{1/2}}\leq C\eps.
\end{align*}
}

\subsection{Results for deterministic and stochastic oracle}

\textbf{Smooth objectives}. Saddle-point problems with deterministic oracle are studied in a more general case (finite-sum structure with variance reduction) in \cite{kovalev2022optimal_1}. The reduction of Algorithm 2 of \cite{kovalev2022optimal_1} for deterministic case achieves optimal complexity bounds for deterministic setup.

Also inexact oracle technique was used in \cite{beznosikov2021distributed_2,beznosikov2021near} to build algorithms optimal up to $\log(1/\eps)$. On the other hand, adding regularization (Lemma~\ref{lem:regularization}) to the results of \cite{kovalev2022optimal_1} allows to avoid squared logarithmic factors and gives an algorithm with {\small$O\cbraces{\chi(\Lmax D^2/\eps)\log(1/\eps)}$} communication and {\small$O\cbraces{(\Lmax D^2/\eps)\log(1/\eps)}$} oracle complexities. In Table~\ref{tab:saddle_deterministic}, \rev{only the results for DESM~\cite{beznosikov2021distributed_2} are given. However, the} regularization of \cite{kovalev2022optimal_1} is also applicable.

\begin{table}[H]
	\begin{tabular}{|b{0.1cm}|b{0.9cm}|>{\centering\arraybackslash}b{2.9cm}|>{\centering\arraybackslash}b{2.9cm}|}
		\hline
		& & convex-concave & str. convex-concave \\
		\hline
		\parbox[t]{20mm}{\multirow{2}{*}{\rotatebox[origin=r]{90}{\footnotesize smooth\hspace{0.2cm}}}} & & DESM \cite{beznosikov2021distributed_2} & Alg.2 \cite{kovalev2022optimal_1} \\
		& ~comm. & $\chi\frac{\Lf D^2}{\eps}\log^2\frac{1}{\eps}$ & $\chi\frac{\Lmax}{\mumin}\log\frac{1}{\eps}$ \\
		& ~comp. & $\frac{\Lf D^2}{\eps}\log\frac{1}{\eps}$ & $\frac{\Lmax}{\mumin}\log\frac{1}{\eps}$ \\
		\hline
		\vspace{0.1cm}
		\parbox[t]{20mm}{\multirow{2}{*}{\rotatebox[origin=c]{90}{\footnotesize{non-smooth}}}} & comm. & \blue{$\chi\frac{\Mav D}{\eps}$} & \blue{$\chi\frac{\Mav}{\sqrt{\mumin\eps}}$} \\ [0.2cm]
		& comp. & \blue{$\frac{\Mav^2D^2}{\eps^2}$} & \blue{$\frac{\Mav^2}{\mumin\eps}$} \\ [0.2cm] \hline
	\end{tabular}
	\caption{Deterministic saddle-point problems}
	\label{tab:saddle_deterministic}
\end{table}

In the stochastic smooth setup, \cite{beznosikov2021distributed_2} proposed a Decentralized extra-step method (DESM) based on inexact oracle framework. Combined with a mini-batching technique, DESM achieves optimal lower bounds in the stochastic case up to a logarithmic factor $\log(1/\eps)$.
\begin{table}[H]
	\begin{tabular}{|b{0.1cm}|b{0.9cm}|>{\centering\arraybackslash}b{2.95cm}|>{\centering\arraybackslash}b{2.9cm}|}
		\hline
		& & convex-concave & str. convex-concave \\ \hline
		\parbox[t]{20mm}{\multirow{2}{*}{\rotatebox[origin=r]{90}{\footnotesize smooth\hspace{0.2cm}}}} & & DESM \cite{beznosikov2021distributed_2} & DESM \cite{beznosikov2021distributed_2} \\
		& ~comm. & $\chi\frac{\Lf D^2}{\eps}\log^2\frac{1}{\eps}$ & $\chi\frac{\Lf}{\muf}\log^2\frac{1}{\eps}$ \\
		& ~comp. & $\max\cbraces{\frac{L_fD^2}{\eps}\log\frac{1}{\eps}, \frac{\sigmamax^2 D^2}{m\eps^2}}$ & $\max\cbraces{\frac{L_f}{\mu_f}\log\frac{1}{\eps}, \frac{\sigmamax^2}{m\muf\eps}}$ \\
		\hline
		\vspace{0.1cm}
		\parbox[t]{20mm}{\multirow{2}{*}{\rotatebox[origin=r]{90}{\footnotesize non-smooth\hspace{-0.1cm}}}} & ~comm. & \blue{$\chi\frac{\Mav D}{\eps}$} & \blue{$\chi\frac{\Mav}{\sqrt{\mumin\eps}}\log\frac{1}{\eps}$} \\ [0.3cm]
		& ~comp. & \blue{$\max\cbraces{\frac{\Mav^2 D^2}{\eps^2}, \frac{\sigmaav^2 D^2}{m\eps^2}}$} & \blue{$\max\cbraces{\frac{\Mav^2}{\mumin\eps}, \frac{\sigmaav^2}{m\mumin\eps}}$} \\ [0.3cm]
		\hline
	\end{tabular}
	\caption{Stochastic saddle-point problems}
	\label{tab:saddle_stoch}
\end{table}
\textbf{Non-smooth objectives}. Non-smooth min-max problems can be solved with a subgradient method \cite{nedic2009subgradient}. However, it is not optimal in communications and therefore nonsmooth distributed saddles are an open direction of research. \rev{It may be possible to} generalize gradient sliding for optimization \cite{dvinskikh2021decentralized} and apply it to saddle-point problems as it was done for time-static graphs in \cite{kuruzov2022gradient}. \rev{The} guesses for deterministic and stochastic objectives are presented in Tables~\ref{tab:saddle_deterministic} and \ref{tab:saddle_stoch}.

\subsection{Finite-sum problems}

Consider saddle-point problems with finite-sum structure of every function held at the node
\begin{align}\label{eq:sum_type_saddle_finite_sum}
	\min_{x\in\R^d} f(x, y) = \frac{1}{m}\sum_{i=1}^m \frac{1}{n}\sum_{j=1}^n f_{ij}(x, y).
\end{align}
Analogically to Section~\ref{subsec:finite_sum_optimization}, let $f_{ij}$ be $L_{ij}$-smooth and define worst-case smoothness constant over summands $\Lsummand$ as in \eqref{eq:def_Lsummand}.

An optimal variance reduction algorithm for saddle-point problems was proposed in \cite{kovalev2022optimal_1}.

\begin{table}[H]
	\begin{tabular}{|b{1.0cm}|>{\centering\arraybackslash}b{3.1cm}|>{\centering\arraybackslash}b{3.1cm}|}
		\hline
		& convex-concave & str.convex-concave \\
		\hline
		& \blue{Alg.2 \cite{kovalev2022optimal_1}} & Alg.2 \cite{kovalev2022optimal_1} \\ 
		~comm. & \blue{$O\cbraces{\chi\frac{\Lsummand R^2}{\eps}\log\frac{1}{\eps}}$} & $O\cbraces{\chi\frac{\Lsummand}{\muav}\log\frac{1}{\eps}}$ \\
		~comp. & \blue{$O\cbraces{\sqrt{n}\frac{\Lsummand R^2}{\eps}\log\frac{1}{\eps}}$} & $O\cbraces{\sqrt{n}\frac{\Lsummand}{\muav}\log\frac{1}{\eps}}$ \\
		\hline
	\end{tabular}
	\caption{Finite-sum smooth min-max problems}
	\label{tab:finite_sum_saddle}
\end{table}

Note that it may be possible to replace $\Lsummand$ with $\Lmax$ by using importance sampling \cite{csiba2018importance}.

\subsection{Different strong convexity and strong concavity constants}

In Tables~\ref{tab:saddle_deterministic} and \ref{tab:saddle_stoch} the results for min-max problems with same constants of strong convexity and strong concavity \rev{were described}. However, strong convexity constant $\mu_x$ may be different from strong concavity constant $\mu_y$.
In the non-distributed case, lower bounds for saddle-point problems with different constants of strong convexity and strong concavity were derived in \cite{zhang2019lower}. Optimal algorithms achieving these bounds were developed for saddles with bilinear coupling, i.e. problems of form
\begin{align*}
	\min_{x\in\R^{d_x}} \max_{y\in\R^{d_y}}~ g(x) + y^\top A x - q(y),
\end{align*}
where $g(x)$ and $q(y)$ are strongly convex functions with constants $\mu_x$ and $\mu_y$, respectively. For saddles with bilinear coupling, optimal algorithms were developed in \cite{jin2022sharper,kovalev2021accelerated,thekumparampil2022lifted}. Paper \cite{carmon2022recapp} proposed a Catalyst type algorithm for general (non-strongly) convex strongly-concave saddles.

Distributed algorithms for saddle point problems with different strong convexity and strong concavity constants are an open venue for research. Paper \cite{metelev2022decentralized} uses an inexact oracle framework and a mini-batching technique and reaches lower bounds up to logarithmic factor.

\rev{\section{Conclusions}

Decentralized optimization has numerous applications in distributed control and sensing, large-scale and privacy preserving machine learning and distributed statistical inference. Developing optimization methods robust to network architecture changes is also an important area of research. This paper surveys the known results in the area and discusses the techniques already known in the literature. However, there is still a number of open questions. In particular, it seems that no universal and simple technique to convert any optimization algorithm to a decentralized algorithm has been studied.
}


\section{Acknowledgments}

The authors are grateful to Angelia Nedi\'c, Gesualdo Scutari and Anton Proskurnikov.

The work of A.~Rogozin and A.~Beznosikov in sections 1, 4 was supported by the strategic academic leadership program <<Priority 2030>> (Agreement  075-02-2021-1316 30.09.2021).

The work of A.~Rogozin and A.~Gasnikov in sections 2, 3 was supported by the Ministry of Science and Higher Education of the Russian Federation (Goszadaniye) 075-00337-20-03, project No. 0714-2020-0005.

\bibliographystyle{abbrv}
\bibliography{references}

\end{document}